\documentclass{article}
\begin{document}
\newtheorem{proposition}{Proposition}[section]
\newtheorem{definition}{Definition}[section]
\newtheorem{lemma}{Lemma}[section]
 
\title{\bf $\xi$-Groups and Hu-Liu Leibniz Algebras}
\author{Keqin Liu\\Department of Mathematics\\The University of British Columbia\\Vancouver, BC\\
Canada, V6T 1Z2}
\date{December 17, 2005}
\maketitle

\begin{abstract} We initiate the study of $\xi$-groups and Hu-Liu Leibniz algebras, claim that almost all simple Leibniz algebras and simple Hu-Liu Leibniz algebras are linear, and establish two passages. One is the passage from a special $\mathcal{Z}_2$-graded associative algebra to a Hu-Liu Leibniz algebra. The other one is the passage from a linear $\xi$-group to its tangent space which is a Hu-Liu Leibniz algebra.
\end{abstract}

The Lie correspondence between connected linear groups and linear Lie algebras is a central result of the Lie theory. Anyone who attempts at generalizing the Lie correspondence is bothered by the following problem: which kinds of algebraic objects should be used to replace groups and Lie algebras? After searching for the generalization of the Lie correspondence for many years, I am still not sure whether it is practical to live in hope of getting a complete generalization of the Lie correspondence by replacing groups and Lie algebras with some kinds of algebraic objects which are not groups and Lie algebras respectively. However, it is certain that the Lie correspondence and the basic Lie theory can be extended completely if some kinds of algebraic structures are added to groups and Lie algebras. In order to extend the Lie correspondence, the simplest replacements of groups and Lie algebras are $\xi$-groups and Hu-Liu Leibniz algebras, which are the algebraic objects obtained by adding more algebraic structures to groups and Lie algebras, respectively. The purpose of this paper is to give the basic properties of $\xi$-groups and Hu-Liu Leibniz algebras.

This paper is organized as follows. First, in Section 1 I review some notions about Leibniz algebras, introduce the notion of a linear Leibniz algebra by using a special $\mathcal{Z}_2$-graded associative algebra, and claim that almost all simple Leibniz algebras are linear. Next, in Section 2 I introduce the notion of a Hu-Liu Leibniz algebra by adding a Lie algebra structure to a Leibniz algebra, establish the passage from a special $\mathcal{Z}_2$-graded associative 
algebra to a Hu-Liu Leibniz algebra, and claim that almost all simple Hu-Liu Leibniz algebras are linear. Finally, in Section 3 I introduce the notion of a $\xi$-group, and establish the passage from a linear $\xi$-group to its tangent space which is a Hu-Liu Leibniz algebra.

\bigskip
\section{Linear Leibniz Algebras}

We begin this section by recalling some concepts about Leibniz algebras in \cite{LT}. 

\medskip
A vector space $L$ over a field $\mathbf{k}$ is called a {\bf right Leibniz algebra} if there exists a binary operation $\langle \, , \, \rangle$: $L\times L\to L$, called the {\bf angle bracket},  such that the {\bf right Leibniz identity} holds:
\begin{equation}\label{eq1}
\langle\langle x, y\rangle , z\rangle=\langle x, \langle y, z\rangle\rangle+ 
\langle\langle x,  z\rangle , y\rangle \qquad\mbox{for $x,y,z\in\mathcal{L}$.}
\end{equation}

Let $L$ be a right Leibniz algebra. A subspace $I$ of $L$ is called an {\bf ideal} of $L$ if 
$\langle I, L\rangle\subseteq I$ and $\langle L, I\rangle\subseteq I$. The {\bf annihilator} $L^{ann}$ of $L$ is defined by
\begin{equation}\label{eq2}
L^{ann}:=\sum _{x\in L}\mathbf{k}\langle x, x\rangle.
\end{equation}
It is easy to check that $L^{ann}$ is an ideal of $L$ and $L^{ann}$ can also be expressed as follows:
\begin{equation}\label{eq3}
L^{ann}:=\sum _{x, y\in L}\mathbf{k}\big(\langle x, y\rangle + \langle y, x\rangle \big).
\end{equation}

A right Leibniz algebra $L$ is said to be {\bf simple} if $L^{ann}\ne 0$ and $L$ has no ideals which are not equal to $\{0\}$, $L^{ann}$ and $L$.

Let $L$ and $\bar{L}$ be two right Leibniz algebras. A linear map $\phi$: $L\to \bar{L}$ is called a {\bf homomorphism} if $\phi\big(\langle x, y\rangle\big)=\langle \phi(x), \phi(y)\rangle$ for $x$, $y\in L$.

\medskip
We now introduce the notion of a special $\mathcal{Z}_2$-graded associative algebra.

\begin{definition}\label{def1.1} An associative algebra $A$ is called  a {\bf special $\mathcal{Z}_2$-graded associative algebra} if $A=A_0\oplus A_1$ (as vector spaces) and 
\begin{equation}\label{eq0}
A_0A_0\subseteq A_0,\quad A_0A_1+A_1A_0\subseteq A_1 \quad\mbox{and}\quad A_1A_1=0.
\end{equation}
\end{definition}

Let $A=A_0\oplus A_1$ be a special $\mathcal{Z}_2$-graded associative algebra, where $A_0$ and $A_1$ satisfy (\ref{eq0}). By Definition~\ref{def1.1}, $A_0$ is a subalgebra of $A$ and $A_1$ is a bimodule over $A_0$. $A_0$ and $A_1$ are called the {\bf even part} and {\bf odd part} of $A$, respectively. An element $a$ of $A$ can be written uniquely as $a=a_0+a_1$, where 
$a_0\in A_0$ and $a_1\in A_1$ are called the {\bf even component} and {\bf odd component} of $a$, respectively.

The next proposition establishes the passage from a special $\mathcal{Z}_2$-graded associative algebra to a Leibniz algebra. 

\begin{proposition}\label{pr1.1} If $A=A_0\oplus A_1$ is a special $\mathcal{Z}_2$-graded associative algebra, then
$( A,\, +,\, \langle \,, \,\rangle)$ is a right Leibniz algebra, where the angle bracket 
$\langle \, , \, \rangle$ is defined by
\begin{equation}\label{eq4}
\langle x, y\rangle: =xy_0-y_0x
\end{equation}
for $x\in A$, $y=y_0+y_1\in A$ and $y_i\in A_i$ with $i=0$ and $1$.
\end{proposition}

\medskip
\noindent
{\bf Proof} It is a straightforward computation.

\hfill\raisebox{1mm}{\framebox[2mm]{}}

\bigskip
\begin{definition}\label{def1.2} A right Leibniz algebra $L$ is said to be {\bf linear} if there exists a special $\mathcal{Z}_2$-graded associative algebra $A$ such that there is an injective homomorphism from the right Leibniz algebra $L$ to the right Leibniz algebra $( A,\, +,\, \langle \,, \,\rangle)$ defined by (\ref{eq4}).
\end{definition}

One important property of a simple right Leibniz algebra is the following

\begin{proposition}\label{pr1.2} If $( L,\, +,\, \langle \,, \,\rangle)$ is a simple right Leibniz algebra such that 
\newline $\langle L^{ann}, L \rangle\ne 0$, then $L$ is linear.
\hfill\raisebox{1mm}{\framebox[2mm]{}}
\end{proposition}

\bigskip
The key of proving Proposition~\ref{pr1.2} is to find the counterpart of the adjoint representation of a Lie algebra in the context of Leibniz algebra. The counterpart of the ordinary adjoint representation has been used to prove the counterparts of Engle's Theorem and Lie's Theorem in the context of Leibniz algebra, which were announced in Appendix 1 of \cite{Liu2}.

\bigskip
\section{Hu-Liu Leibniz Algebras}

We now introduce the abstract notion of a Hu-Liu Leibniz algebra in the following

\begin{definition}\label{def2.1} A right Leibniz algebra $( L,\, +,\, \langle \,, \,\rangle)$ is called a {\bf right Hu-Liu Leibniz algebra} if there exists a binary operation 
$[\, , \, ] : L\times L \to L$ such that $( L,\, +,\, [\, ,\, ])$ is a Lie algebra and the following {\bf right Hu-Liu identities} hold:
\begin{equation}\label{eq5}
\langle x, [y, z]\rangle =\langle x, \langle y, z\rangle\rangle,
\end{equation}
\begin{equation}\label{eq6}
[\langle x, x\rangle , y] =\langle \langle x, x\rangle , y\rangle,
\end{equation}
\begin{equation}\label{eq7}
\langle [x, y], z\rangle +[\langle y, z\rangle , x]+[y, \langle x, z\rangle ]=0,
\end{equation}
\begin{equation}\label{eq8}
[\langle x, y\rangle , z] +[z, [x, y]]+[z, \langle y, x\rangle  ]+
\langle z, \langle x, y\rangle\rangle=0,
\end{equation}
where $x$, $y$, $z\in L$.
\end{definition}

The notion of a left Hu-Liu Leibniz algebra can be introduced similarly. In this paper, we deal with only right Hu-Liu Leibniz algebras. From now on, a Hu-Liu Leibniz algebra $L$ always means a right Hu-Liu Leibniz algebra and is denoted by $( L,\, +,\, \langle \,, \,\rangle ,\, [\, ,\, ])$, where 
$\langle \,, \,\rangle $ and $ [\, ,\, ]$ are called the {\bf angle bracket} and the {\bf square bracket}, respectively.

\begin{definition}\label{def2.2} Let $I$ be a subspace of a Hu-Liu Leibniz algebra 
$( L,\, +,\, \langle \,, \,\rangle ,\, [\, ,\, ])$.
\begin{description}
\item[(i)] $I$ is called an {\bf ideal} of $L$ if
$$
\langle I, L\rangle\subseteq I, \quad \langle L, I\rangle\subseteq I \quad\mbox{and}\quad 
[L, I]\subseteq I.
$$
\item[(ii)] $I$ is called a {\bf Hu-Liu Leibniz subalgebra} of $L$ if
$$
\langle I, I\rangle\subseteq I \quad\mbox{and}\quad [I, I]\subseteq I.
$$
\end{description}
\end{definition}

A Hu-Liu Leibniz algebra $L$ always has at least three ideals: $\{0\}$, $L^{ann}$ and $L$ by (\ref{eq6}). 

\begin{definition}\label{def2.3} A Hu-Liu Leibniz algebra $L$ is called a {\bf simple Hu-Liu Leibniz algebra} if $L^{ann}\ne 0$ and $L$ has no ideals which are not equal to $\{0\}$, $L^{ann}$ and $L$.
\end{definition}

If $( L,\, +,\, \langle \,, \,\rangle ,\, [\, ,\, ])$ is a Hu-Liu Leibniz algebra, then $[L^{ann}, L^{ann}]=0$; that is, $L^{ann}$ is an Abelian ideal of the Lie algebra $( L,\, +,\, [\, ,\, ])$. Hence, if $L^{ann}\ne 0$, then the Lie algebra $( L,\, +,\, [\, ,\, ])$ carried by a Hu-Liu Leibniz algebra $L$ is always nonsemisimple.

\begin{definition}\label{def2.4} Let $L$ and $\bar{L}$ be two Hu-Liu Leibniz algebras. A linear map $\phi$: $L\to \bar{L}$ is called a {\bf homomorphism} if 
$$\phi\big(\langle x, y\rangle\big)=\langle \phi(x), \phi(y)\rangle \quad\mbox{and}\quad 
\phi\big([ x, y]\big)=[ \phi(x), \phi(y)] $$
for $x$, $y\in L$.
\end{definition}

Clearly, if $\phi$: $L\to \bar{L}$ is a Hu-Liu Leibniz algebra homomorphism, then the {\bf kernel}
$$Ker\phi :=\{\, x\in L \, |\, \phi (x)=0 \,\}$$
is an ideal of $L$, and the {\bf image}
$$ Im\phi : =\{\, \phi (x) \, |\, x\in L \, \}$$
is a Hu-Liu Leibniz subalgebra of $\bar{L}$.

The next proposition establishes the passage from a special $\mathcal{Z}_2$-graded associative algebra to a Hu-Liu Leibniz algebra. 

\begin{proposition}\label{pr2.1} If $A=A_0\oplus A_1$ is a special $\mathcal{Z}_2$-graded associative algebra, then
$( A,\, +,\, \langle \,, \,\rangle , $ $\, [\, ,\, ])$ is a Hu-Liu Leibniz algebra, where the angle bracket $\langle \, , \, \rangle$ and the square bracket $[\, , \,]$ are defined by
\begin{equation}\label{eq9}
\langle x, y\rangle: =xy_0-y_0x \quad\mbox{and}\quad [x, y]:=xy-yx
\end{equation}
for $x\in A$, $y=y_0+y_1\in A$ and $y_i\in A_i$ with $i=0$ and $1$.
\end{proposition}

\medskip
\noindent
{\bf Proof} It is a straightforward computation.

\hfill\raisebox{1mm}{\framebox[2mm]{}}

\bigskip
The angle bracket $\langle \, , \, \rangle$ and the square bracket $[\, , \,]$ defined by  (\ref{eq9}) satisfy a number of identities. The identities (\ref{eq5}), (\ref{eq6}), (\ref{eq7}) and (\ref{eq8}) in Definition~\ref{def2.1} are just four of them. 

\begin{definition}\label{def2.5} A Hu-Liu Leibniz algebra $L$ is said to be {\bf linear} if there exists a special $\mathcal{Z}_2$-graded associative algebra $A$ such that there is an injective homomorphism from the Hu-Liu Leibniz algebra $L$ to the Hu-Liu Leibniz algebra 
$( A,\, +,\, \langle \,, \,\rangle , \,[\, , \,])$ defined by (\ref{eq9}).
\end{definition}

Using the same idea of proving Proposition~\ref{pr1.2}, we have

\begin{proposition}\label{pr2.2} If $( L,\, +,\, \langle \,, \,\rangle, \,[\, , \,])$ is a simple Hu-Liu Leibniz algebra such that $[ L^{ann}, L]\ne 0$, then $L$ is linear.
\hfill\raisebox{1mm}{\framebox[2mm]{}}
\end{proposition}

\bigskip
\section{$\xi$-Groups}

The notion of a $\xi$-group is based on the following

\begin{definition}\label{def3.1} A pair $(\mathcal{A}, \xi)$  is called a {\bf covering pair} if $\mathcal{A}$ is a group and $\xi$: $\mathcal{A}\to \mathcal{A}$ is a group homomorphism with $\xi ^2=\xi$.
\end{definition}

We now introduce $\xi$-groups which are the class of groups used in my extension of the Lie correspondence.

\begin{definition}\label{def3.2} A group $G$ is called a {\bf $\xi$-group} if there exists a covering pair $(\mathcal{A}, \xi)$ such that $G$ is a subgroup of $\mathcal{A}$ and
\begin{equation}\label{eq10}
\xi(x) \,G\, \xi(x)^{-1}\subseteq G \quad\mbox{for $x\in G$.}
\end{equation}
\end{definition}

A $\xi$-group $G$ is also denoted by $(G; \mathcal{A}, \xi)$, where the pair $(\mathcal{A}, \xi)$ is also called the {\bf covering pair} of $G$ and $G$ is also called a $\xi$-group in the covering pair $(\mathcal{A}, \xi)$.
 
\medskip
Let $A=A_0\oplus A_1$ be a special $\mathcal{Z}_2$-graded associative algebra with the identity $1$. We define 
$A^{-1}$ and $A_0^{-1}$ by 
$$
A^{-1}:=\{\, x\in A \, |\, \mbox{there exists $y\in A$ such that $yx=1=xy$ }\,\}
$$
and
$$
A_0^{-1}:=\{\, x_0\in A_0 \, |\, \mbox{there exists $y_0\in A_0$ such that $y_0x_0=1=x_0y_0$} \,\}.
$$
It is clear that
$$
A^{-1}=\{\, x_0+x_1 \, |\, \mbox{$x_0\in A_0^{-1}$ and $x_1\in A_1$ }\,\}
$$
and the map
\begin{equation}\label{eq12}
\xi _{A^{-1}}: x_0+x_1 \to x_0 \quad\mbox{for $x_0\in A_0^{-1}$} 
\end{equation}
is a group homomorphism from $A^{-1}$ to $A^{-1}$ with $\xi _{A^{-1}}^2=\xi _{A^{-1}}$. Hence, 
$(A^{-1}, \xi _{A^{-1}})$ is a covering pair. By (\ref{eq10}) and (\ref{eq12}), a subgroup $G$ of 
$A^{-1}$ is a $\xi$-group in the covering pair $(A^{-1}, \xi _{A^{-1}})$ if and only if
$$
x_0\, G \, x_0^{-1}\subseteq G \quad
\mbox{for $x_0+x_1\in G$ with $x_0\in A_0$ and $x_1\in A_1$ .}
$$
A $\xi$-group in the covering pair $(A^{-1}, \xi _{A^{-1}})$ is called a {\bf linear $\xi$-group} in the  special $\mathcal{Z}_2$-graded associative algebra $A$.

\medskip
In the remaining part of this paper, $\mathbf{k}$ will denote the field $\mathcal{R}$ of real numbers or the field of complex numbers.

\begin{definition}\label{def3.3} An associative algebra $A$ over $\mathbf{k}$ is called a 
{\bf special $\mathcal{Z}_2$-graded Banach algebra} if $A$ is both a special $\mathcal{Z}_2$-graded associative algebra and a Banach algebra.
\end{definition}

Let $G$ be a linear $\xi$-group in a special $\mathcal{Z}_2$-graded Banach algebra $A$. The {\bf tangent space} $T_1(G)$ to $G$ at the identity $1$ of $A$ is defined by
$$
T_1(G):=\{\, a'(0) \, |\, \mbox{$a(t): \mathcal{R}\to G$ is a differential curve with $a(0)=1$ }\,\}.
$$ 

\begin{proposition}\label{pr3.1} Let $A$ be a special $\mathcal{Z}_2$-graded Banach algebra over a field $\mathbf{k}$. If $G$ is a linear $\xi$-group in $A$, then the tangent space $T_1(G)$ to $G$ at the identity $1$ of $A$ is a Hu-Liu Leibniz algebra over the field of real numbers.
\end{proposition}

\medskip
\noindent
{\bf Proof} The proof is the same as the proof of Proposition 4.4 in \cite{Liu2}. 

\hfill\raisebox{1mm}{\framebox[2mm]{}}

\bigskip
The tangent space $T_1(G)$ to a linear $\xi$-group $G$ is also called the {\bf Hu-Liu Leibniz algebra of $G$}. Proposition~\ref{pr3.1} establishes the passage from a linear $\xi$-group $G$ to its Hu-Liu Leibniz algebra. 

\bigskip
I end the paper with two remarks. One is that the passage in
Proposition~\ref{pr3.1} is in fact an one-to-one correspondence between connected linear $\xi$-groups and linear Hu-Liu Leibniz algebras.
The other one is that Frobenius' theorem about the classification of finite dimensional division real associative algebras has a natural extension, which is an unexpected result appearing in my search for the generalization of the Lie correspondence. It turns out that if the notion of a division associative algebra is modified in the context of special $\mathcal{Z}_2$-graded associative algebras, then the classification of finite dimensional modified division real associative algebras consists of eight objects, where the first three objects are the ordinary finite dimensional division real associative algebras (without divisors of zero), and the last five objects are finite dimensional special $\mathcal{Z}_2$-graded associative algebras with divisors of zero. These results will be announced in \cite{hl2}.

\end{document}